\documentclass{amsart}
\pdfoutput=1

\usepackage{amsmath,amssymb}
\usepackage{microtype}
\usepackage{booktabs}
\usepackage{url}

\newtheorem{theorem}{Theorem}

\numberwithin{equation}{section}

\title{Approximate common divisors via lattices}

\author{Henry Cohn}
\address{Microsoft Research New England\\
One Memorial Drive\\
Cambridge, MA 02142}
\email{cohn@microsoft.com}

\author{Nadia Heninger}
\address{Department of Computer Science\\
Princeton University\\
Princeton, NJ 08540}
\curraddr{Department of Computer Science and Engineering\\
University of California, San Diego\\
9500 Gilman Drive, Mail Code 0404\\
La Jolla, CA 92093-0404}
\email{nadiah@cs.ucsd.edu}


\date{March 13, 2012}

\begin{document}

\begin{abstract}
  We analyze the multivariate generalization of Howgrave-Graham's
  algorithm for the approximate common divisor problem.  In the
  $m$-variable case with modulus $N$ and approximate common divisor of
  size $N^\beta$, this improves the size of the error tolerated from
  $N^{\beta^2}$ to $N^{\beta^{(m+1)/m}}$, under a commonly used
  heuristic assumption.  This gives a more detailed analysis of the
  hardness assumption underlying the recent fully homomorphic
  cryptosystem of van Dijk, Gentry, Halevi, and Vaikuntanathan.  While
  these results do not challenge the suggested parameters, a
  $2^{n^{\varepsilon}}$ approximation algorithm with $\varepsilon<2/3$
  for lattice basis reduction in $n$ dimensions could be used to break
  these parameters.  We have implemented our algorithm, and it
  performs better in practice than the theoretical analysis suggests.

  Our results fit into a broader context of analogies between
  cryptanalysis and coding theory.
  The multivariate approximate common divisor problem is the
  number-theoretic analogue of multivariate polynomial
  reconstruction, and we develop a corresponding lattice-based
  algorithm for the latter problem.  In particular, it specializes to
  a lattice-based list decoding algorithm for Parvaresh-Vardy and
  Guruswami-Rudra codes, which are multivariate extensions of
  Reed-Solomon codes.  This yields a new proof of the list decoding
  radii for these codes.
\end{abstract}

\maketitle

\section{Introduction}

Given two integers, we can compute their greatest common divisor
efficiently using Euclid's algorithm.  Howgrave-Graham
\cite{howgrave2001approximate} formulated and gave an algorithm to
solve an approximate version of this question, asking the question
``What if instead of exact multiples of some common divisor, we only
know approximations?''  In the simplest case, we are given one exact
multiple $N = p q_0$ and one near multiple $a_1 = p q_1 + r_1$, and the
goal is to learn $p$, or at least $p \gcd(q_0,q_1)$.

In this paper, we generalize Howgrave-Graham's approach to the case
when one is given many near multiples of $p$.  The hardness of solving
this problem for small $p$ (relative to the size of the near multiples)
was recently proposed as the foundation for a fully homomorphic
cryptosystem \cite{vandijkgentryhalevivaikuntanathan2010}.
Specifically, we can show that improving the approximation of lattice
basis reduction for the particular lattices $L$ we are looking at from
$2^{\dim L}$ to $2^{(\dim L)^{\varepsilon}}$ with $\varepsilon < 2/3$
would break the suggested parameters in the system. See
Section~\ref{section:fail} for the details.  The approximate common
divisor problem is also closely related to the problem of finding small
solutions to multivariate polynomials, a problem first posed by
Coppersmith~\cite{coppersmith1997small}, and whose various extensions
have many applications in cryptanalysis~\cite{boneh2000cryptanalysis}.

The multivariate version of the problem allows us to improve the bounds
for when the approximate common divisor problem is solvable: given
$N=pq_0$ and $m$ randomly chosen approximate multiples $a_i = p q_i +
r_i$ of $p = N^\beta$, as well as upper bounds $X_i$ for each $|r_i|$,
we can find the perturbations $r_i$ when
\[
\sqrt[m]{X_1 \dots X_m}
< N^{(1+o(1))\beta^{(m+1)/m}}.
\]
In other words, we can compute approximate common divisors when $r_i$
is as large as $N^{\beta^{(m+1)/m}}$.  For $m=1$, we recover
Howgrave-Graham's theorem \cite{howgrave2001approximate}, which handles
errors as large as $N^{\beta^2}$.  As the number $m$ of samples grows
large, our bound approaches $N^\beta$, i.e., the size of the
approximate common divisor $p$.  Our algorithm runs in polynomial time
for fixed $m$.  We cannot rigorously prove that it always works, but it
is supported by a heuristic argument and works in practice.

There is an analogy between the ring of integers and the ring of
polynomials over a field.  Under this analogy, finding a large
approximate common divisor of two integers is analogous to
reconstructing a polynomial from noisy interpolation information, as we
explain in Section~\ref{section:multipol}.  One of the most important
applications of polynomial reconstruction is decoding of Reed-Solomon
codes.  Guruswami and Sudan \cite{guruswami1998improved} increased the
feasible decoding radius of these codes by giving a list-decoding
algorithm that outputs a list of polynomially many solutions to a
polynomial reconstruction problem.  The analogy between the integers
and polynomials was used in \cite{ideal-coppersmith} to give a proof of
the Guruswami-Sudan algorithm inspired by Howgrave-Graham's approach,
as well as a faster algorithm.

Parvaresh and Vardy \cite{parvaresh2005correcting} developed a related
family of codes with a larger list-decoding radius than Reed-Solomon
codes.  The decoding algorithm corresponds to simultaneous
reconstruction of several polynomials.

In this paper, we observe that the problem of simultaneous
reconstruction of multiple polynomials is the exact analogue of the
approximate common divisor problem with many inputs, and the improved
list-decoding radius of Parvaresh-Vardy codes corresponds to the
improved error tolerance in the integer case.  We adapt our algorithm
for the integers to give a corresponding algorithm to solve the
multiple polynomial reconstruction problem.

This algorithm has recently been applied to construct an optimally
Byzantine-robust private information retrieval protocol~\cite{DGH}. The
polynomial lattice methods we describe are extremely fast in practice,
and they speed up the client-side calculations by a factor of several
thousand compared with a related scheme that uses the Guruswami-Sudan
algorithm.  See \cite{DGH} for more information and timings.

\subsection{Related work}

Howgrave-Graham first posed the problem of approximate integer common
divisors in \cite{howgrave2001approximate}, and used it to address the
problem of factoring when information is known about one of the
factors.  His algorithm gave a different viewpoint on Coppersmith's
proof \cite{coppersmith1997small} that one can factor an RSA modulus $N
= pq$ where $p \approx q \approx \sqrt{N}$ given the most significant
half of the bits of one of the factors.  This technique was applied by
Boneh, Durfee, and Howgrave-Graham \cite{boneh1999factoring} to factor
numbers of the form $p^r q$ with $r$ large. Jochemsz and
May~\cite{jochemszmay} and Jutla~\cite{jutla1998finding} considered the
problem of finding small solutions to multivariate polynomial
equations, and showed how to do so by obtaining several equations
satisfied by the desired roots using lattice basis reduction. Herrmann
and May \cite{herrmann2008solving} gave a similar algorithm in the case
of finding solutions to multivariate linear equations modulo divisors
of a given integer. They applied their results to the case of factoring
with bits known when those bits might be spread across $\log\log N$
chunks of $p$.  Notably, their results display similar behavior to ours
as the number of variables grows large.  Van Dijk, Gentry, Halevi, and
Vaikuntanathan \cite{vandijkgentryhalevivaikuntanathan2010} discuss
extensions of Howgrave-Graham's method to larger $m$ and provide a
rough heuristic analysis in Appendix~B.2 of the longer version of their
paper available on the Cryptology ePrint Archive.

Chen and Nguyen~\cite{chen2011faster} gave an algorithm to find
approximate common divisors which is not related to the
Coppersmith/Howgrave-Graham lattice techniques and which provides an
exponential speedup over exhaustive search over the possible
perturbations.

In addition to the extensive work on polynomial reconstruction and
noisy polynomial interpolation in the coding theory literature, the
problem in both the single and multiple polynomial cases has been used
as a cryptographic primitive, for example in \cite{kiayias2001secure},
\cite{kiayias2001polynomial}, and \cite{augot2003public} (broken in
\cite{coron2004cryptanalysis}).  Coppersmith and Sudan
\cite{coppersmith2003reconstructing} gave an algorithm for simultaneous
reconstruction of multiple polynomials, assuming random (rather than
adversarially chosen) errors.  Bleichenbacher, Kiayias, and
Yung~\cite{bleichenbacher2003noisy} gave a different algorithm for
simultaneous reconstruction of multiple polynomials under a similar
probabilistic model.  Parvaresh and Vardy
\cite{parvaresh2005correcting} were the first to beat the list-decoding
performance of Reed-Solomon codes for adversarial errors, by combining
multiple polynomial reconstruction with carefully chosen constraints on
the polynomial solutions; this allowed them to prove that their
algorithm ran in polynomial time, without requiring any heuristic
assumptions.  Finally, Guruswami and Rudra \cite{guruswami2006explicit}
combined the idea of multi-polynomial reconstruction with an optimal
choice of polynomials to construct codes that can be list-decoded up to
the information-theoretic bound (for large alphabets).

\subsection{Problems and results}

\subsubsection{Approximate common divisors}
\label{subsub:acd}

Following Howgrave-Graham, we define the ``partial'' approximate common
divisor problem to be the case when one has $N=pq_0$ and $m$
approximate multiples $a_i = pq_i+r_i$ of $p$.  We want to recover an
approximate common divisor.  To do so, we will compute $r_1,\dots,r_m$,
after which we can simply compute the exact greatest common divisor of
$N,a_1-r_1,\dots,a_m-r_m$.

If the perturbations $r_i$ are allowed to be as large as $p$, then it
is clearly impossible to reconstruct $p$ from this data.  If they are
sufficiently small, then one can easily find them by a brute force
search.  The following theorem interpolates between these extremes: as
$m$ grows, the bound on the size of $r_i$ approaches the trivial upper
bound of $p$.

\begin{theorem}[Partial approximate common divisors]
\label{pseudo:pacd} Given positive integers $N,a_1,\dots,a_m$ and
bounds $\beta \gg 1/\sqrt{\log N}$ and $X_1,\dots,X_m$, we can find all
$r_1,\dots,r_m$ such that
\[
\gcd(N,a_1-r_1,\dots,a_m-r_m) \ge N^\beta
\]
and $|r_i| \le X_i$, provided that
\[
\sqrt[m]{X_1 \dots X_m} < N^{(1+o(1))\beta^{(m+1)/m}}
\]
and that the algebraic independence hypothesis discussed in
Section~\ref{section:acd} holds.  The algorithm runs in polynomial time
for fixed $m$, and the $\gg$ and $o(1)$ are as $N \to \infty$.
\end{theorem}

For $m=1$, this theorem requires no algebraic independence hypothesis
and is due to Howgrave-Graham \cite{howgrave2001approximate}. For
$m>1$, not all inputs $N,a_1,\dots,a_m$ will satisfy the hypothesis.
Specifically, we must rule out attempting to improve on the $m=1$ case
by deriving $a_2,\dots,a_m$ from $a_1$, for example by taking $a_i$ to
be a small multiple of $a_1$ plus an additional perturbation (or, worse
yet, $a_1=\dots=a_m$).  However, we believe that generic integers will
work, for example integers chosen at random from a large range, or at
least integers giving independent information in some sense.

We describe our algorithm to solve this problem in
Section~\ref{section:acd}. We follow the general technique of
Howgrave-Graham: we use LLL lattice basis reduction to construct $m$
polynomials for which $r_1,\dots,r_m$ are roots, and then we solve the
system of equations. The lattice basis reduction is for a lattice of
dimension at most $\beta \log N$, regardless of what $m$ is, but the
root finding becomes difficult when $m$ is large.

This algorithm is heuristic, because we assume we can obtain $m$ short
lattice vectors representing algebraically independent polynomials from
the lattice that we will construct.  This assumption is commonly made
when applying multivariate versions of Coppersmith's method, and has
generally been observed to hold in practice.  See
Section~\ref{section:acd} for more details.  This is where the
restriction to generic inputs becomes necessary: if $a_1,\dots,a_m$ are
related in trivial ways, then the algorithm will simply recover the
corresponding relations between $r_1,\dots,r_m$, without providing
enough information to solve for them.

Note that we are always able to find one nontrivial algebraic relation
between $r_1,\dots,r_m$, because LLL will always produce at least one
short vector.  If we were provided in advance with $m-1$ additional
relations, carefully chosen to ensure that they would be algebraically
independent of the new one, then we would have no need for heuristic
assumptions.  We will see later in this section that this situation
arises naturally in coding theory, namely in Parvaresh-Vardy codes
\cite{parvaresh2005correcting}.

The condition $\beta \gg 1/\sqrt{\log N}$ arises from the exponential
approximation factor in LLL.  It amounts to $N^{\beta^2} \gg 1$.  An
equivalent formulation is $\log p \gg \sqrt{\log N}$; i.e., the number
of digits in the approximate common factor $p$ must be more than the
square root of the number of digits in $N$.  When $m=1$, this is not a
restriction at all: when $p$ is small enough that $N^{\beta^2}$ is
bounded, there are only a bounded number of possibilities for $r_1$ and
we can simply try all of them.  When $m>1$, the multivariate algorithm
can handle much larger values of $r_i$ for a given $p$, but the $\log p
\gg \sqrt{\log N}$ condition dictates that $p$ cannot be any smaller
than when $m=1$. Given a lattice basis reduction algorithm with
approximation factor $2^{(\dim L)^{\varepsilon}}$, one could replace
this condition with $\beta^{1+\varepsilon} \log N \gg 1$.  If
$\varepsilon=1/m$, then the constraint could be removed entirely in the
$m$-variable algorithm.  See Section~\ref{section:acd} for the details.

The $\log p \gg \sqrt{\log N}$ condition is the only thing keeping us
from breaking the fully homomorphic encryption scheme from
\cite{vandijkgentryhalevivaikuntanathan2010}.  Specifically, improving
the approximation of lattice basis reduction for the particular
lattices $L$ we are looking at to $2^{(\dim L)^{\varepsilon}}$ with
$\varepsilon < 2/3$ would break the suggested parameters in the system.
See Section~\ref{section:fail} for the details.

We get nearly the same bounds for the ``general'' approximate common
divisor problem, in which we are not given the exact multiple $N$.

\begin{theorem}[General approximate common divisors]
\label{pseudo:gacd} Given positive integers $a_1,\dots,a_m$ (with $a_i
\approx N$ for all $i$) and bounds $\beta \gg 1/\sqrt{\log N}$ and $X$,
we can find all $r_1,\dots,r_m$ such that
\[
\gcd(a_1-r_1,\dots,a_m-r_m) \ge N^\beta
\]
and $|r_i| \le X$, provided that
\[
X< N^{(C_m+o(1))\beta^{m/(m-1)}},
\]
where
\[
C_m = \frac{1-1/m^2}{m^{1/(m-1)}} \approx 1 - \frac{\log m}{m},
\]
and that the algebraic independence hypothesis holds. The algorithm
runs in polynomial time for fixed $m$, and the $\gg$ and $o(1)$ are as
$N \to \infty$.
\end{theorem}

Again, for $m=2$, this result is due to Howgrave-Graham
\cite{howgrave2001approximate}, and no algebraic independence
hypothesis is needed.

The proof is very similar to the case when $N$ is known, but the
calculations are more tedious because the determinant of the lattice is
more difficult to bound.  See Section~\ref{section:gacd} for the
details.

In \cite{howgrave2001approximate}, Howgrave-Graham gives a more
detailed analysis of the behavior for $m=2$.  Instead of our exponent
$C_2\beta^2 = \frac{3}{8}\beta^2$, he gets
$1-\beta/2-\sqrt{1-\beta-\beta^2/2}$, which is asymptotic to
$\frac{3}{8}\beta^2$ for small $\beta$ but is slightly better when
$\beta$ is large.  We are interested primarily in the case when $\beta$
is small, so we have opted for simplicity, but one could carry out a
similar analysis for all $m$.

\subsubsection{Noisy multi-polynomial reconstruction}
\label{section:multipol}

Let $F$ be a field.  Given $m$ single-variable polynomials
$g_1(z),\dots,g_m(z)$ over $F$ and $n$ distinct points $z_1,\dots,z_n$
in $F$, evaluating the polynomials at these points yields $mn$ elements
$y_{ij} = g_i(z_j)$ of $F$.

The noisy multi-polynomial reconstruction problem asks for the recovery
of $g_1,\dots,g_m$ given the evaluation points $z_1,\dots,z_n$, degree
bounds $\ell_i$ on $g_i$, and possibly incorrect values $y_{ij}$.
Stated more precisely: we wish to find all $m$-tuples of polynomials
$(g_1, \dots, g_m)$ satisfying $\deg g_i \le \ell_i$, for which there
are at least $\beta n$ values of $j$ such that $g_i(z_j) = y_{ij}$ for
all $i$.  In other words, some of the data may have been corrupted, but
we are guaranteed that there are at least $\beta n$ points at which all
the values are correct.

Bleichenbacher and Nguyen~\cite{bleichenbacher2000noisy} distinguish
the problem of ``polynomial reconstruction'' from the ``noisy
polynomial interpolation'' problem.  Their definition of ``noisy
polynomial interpolation'' involves reconstructing a single polynomial
when there are several possibilities for each value.  The multivariate
version of this problem can be solved using
Theorem~\ref{pseudo:highdeg}.

This problem is an important stepping stone between single-variable
interpolation problems and full multivariate interpolation, in which we
reconstruct polynomials of many variables.  The multi-polynomial
reconstruction problem allows us to take advantage of multivariate
techniques to prove much stronger bounds, without having to worry about
issues such as whether our evaluation points are in general position.

We can restate the multi-polynomial reconstruction problem slightly to
make the analogy with the integer case clear.  Given evaluation points
$z_j$ and values $y_{ij}$, define $N(z) = \prod_j (z-z_j)$, and use
ordinary interpolation to find polynomials $f_i(z)$ such that $f_i(z_j)
= y_{ij}$.  Then we will see shortly that $g_1,\dots,g_m$ solve the
noisy multi-polynomial reconstruction problem iff
\[
\deg \gcd
(f_1(z)-g_1(z), \dots,f_m(z) - g_m(z), N(z)) \geq \beta n.
\]
This is completely analogous to the approximate common divisor problem,
with $N(z)$ as the exact multiple and $f_1(z),\dots,f_m(z)$ as the
approximate multiples.

To see why this works, observe that the equation $g_i(z_j) = y_{ij}$ is
equivalent to $g_i(z) \equiv y_{ij} \pmod{z-z_j}$.  Thus,
$g_i(z_j)=f_i(z_j)=y_{ij}$ iff $f_i(z) - g_i(z) \equiv 0 \pmod{z-z_j}$,
and $\deg \gcd(f_i(z)-g_i(z),N(z))$ counts how many $j$ satisfy
$g_i(z_j) = y_{ij}$.  Finally, to count the $j$ such that $g_i(z_j) =
y_{ij}$ for all $i$, we use
\[
\deg \gcd
(f_1(z)-g_1(z), \dots,f_m(z) - g_m(z), N(z)).
\]

This leads us to our result in the polynomial case.

\begin{theorem}
\label{multi-polynomial}
Given polynomials $N(z),f_1(z),\dots,f_m(z)$ and degree bounds
$\ell_1,\dots,\ell_m$, we can find all $g_1(z),\dots,g_m(z)$ such that
\[
\deg \gcd (f_1(z) - g_1(z),\dots,f_m(z) - g_m(z),N(z)) \geq \beta \deg N(z)
\]
and $\deg g_i \le \ell_i$, provided that
\[
\frac{\ell_1 + \dots + \ell_m}{m} < \beta^{(m+1)/m} \deg N(z)
\]
and that the algebraic independence hypothesis holds.  The algorithm
runs in polynomial time for fixed $m$.
\end{theorem}


As in the integer case, our analysis depends on an algebraic
independence hypothesis, but it may be easier to resolve this issue in
the polynomial case, because lattice basis reduction is far more
effective and easier to analyze over polynomial rings than it is over
the integers.

Parvaresh-Vardy codes~\cite{parvaresh2005correcting} are based on noisy
multi-polynomial reconstruction: a codeword is constructed by
evaluating polynomials $f_1,\dots,f_m$ at points $z_1,\dots,z_n$ to
obtain $mn$ elements $f_i(z_j)$.  In their construction, $f_1, \dots,
f_m$ are chosen to satisfy $m-1$ polynomial relations, so that they
only need to find one more algebraically independent relation to solve
the decoding problem.  Furthermore, the $m-1$ relations are constructed
so that they must be algebraically independent from the relation
constructed by the decoding algorithm.  This avoids the need for the
heuristic assumption discussed above in the integer case. Furthermore,
the Guruswami-Rudra codes~\cite{guruswami2006explicit} achieve improved
rates by constructing a system of polynomials so that only $n$ symbols
need to be transmitted, rather than $mn$.

Parvaresh and Vardy gave a list-decoding algorithm using the method of
Guruswami and Sudan, which constructs a polynomial by solving a system
of equations to determine the coefficients.  In our terms, they proved
the following theorem:

\begin{theorem}
\label{multi-polynomial-rigorous}
Given a polynomial $N(z)$ and $m$ polynomials $f_1(z),\dots,f_m(z)$,
and degree bounds $\ell_1,\dots,\ell_m$,
we can find a nontrivial polynomial $Q(x_1,\dots,x_m)$ with the
following property: for all $g_1(z),\dots,g_m(z)$ such that
\[
\deg \gcd (f_1(z) - g_1(z),\dots,f_m(z) - g_m(z),N(z)) \geq \beta \deg N(z)
\]
and $\deg g_i \le \ell_i$,
we have
\[
Q(g_1(z),\dots,g_m(z)) = 0,
\]
provided that
\[
\frac{\ell_1 + \dots + \ell_m}{m} < \beta^{(m+1)/m} \deg N(z).
\]
The algorithm runs in polynomial time.  
\end{theorem}

In Section~\ref{section:parvaresh-vardy}, we give an alternative proof
of this theorem using the analogue of lattice basis reduction over
polynomial rings.  This algorithm requires neither heuristic
assumptions nor conditions on $\beta$.

\section{Computing approximate common divisors}
\label{section:acd}

In this section, we describe our algorithm to solve the approximate
common divisor problem over the integers.

To derive Theorem~\ref{pseudo:pacd}, we will use the following
approach:
\begin{enumerate}
\item Construct polynomials $Q_1,\dots,Q_m$ of $m$ variables such that
  \[Q_i(r_1, \dots, r_m) = 0\] for all
  $r_1,\dots,r_m$ satisfying the
  conditions of the theorem.
\item Solve this system of equations to learn candidates for the roots
  $r_1,\dots,r_m$.
\item Test each of the polynomially many candidates to see if it is a
  solution to the original problem.
\end{enumerate}

In the first step, we will construct polynomials $Q$ satisfying
\[
Q(r_1, \dots, r_m) \equiv 0 \pmod{p^k}
\]
(for a $k$ to be chosen later) whenever $a_i \equiv r_i \pmod{p}$ for
all $i$.  We will furthermore arrange that
\[
|Q(r_1,\dots,r_m)| < N^{\beta k}.
\]
These two facts together imply that $Q(r_1,\dots,r_m) = 0$ whenever $p
\ge N^{\beta}$.

To ensure that $Q(r_1, \dots, r_m) \equiv 0 \pmod{p^k}$, we will
construct $Q$ as an integer linear combination of products
\[
(x_1-a_1)^{i_1} \dots (x_m - a_m)^{i_m} N^\ell
\]
with $i_1 + \dots + i_m + \ell \ge k$.  Alternatively, we can think of
$Q$ as being in the integer lattice generated by the coefficient
vectors of these polynomials.  To ensure that $|Q(r_1,\dots,r_m)| <
N^{\beta k}$, we will construct $Q$ to have small coefficients; i.e.,
it will be a short vector in the lattice.

More precisely, we will use the lattice $L$ generated by the coefficient
vectors of the polynomials
\[
(X_1x_1-a_1)^{i_1} \dots (X_mx_m-a_m)^{i_m} N^\ell
\]
with $i_1 + \dots + i_m \le t$ and $\ell = \max\big(k - \sum_j
i_j,0\big)$.  Here $t$ and $k$ are parameters to be chosen later.
Note that we have incorporated the bounds $X_1,\dots,X_m$ on the
desired roots $r_1,\dots,r_m$ into the lattice.  We define $Q$ to be
the corresponding integer linear combination of $(x_1-a_1)^{i_1} \dots
(x_m-a_m)^{i_m} N^\ell $, without $X_1,\dots,X_m$.

Given a polynomial $Q(x_1,\dots,x_m)$ corresponding to a vector $v \in
L$, we can bound
$|Q(r_1,\dots,r_m)|$ by the $\ell_1$ norm $|v|_1$.  Specifically,
if
\[
Q(x_1,\dots,x_m) = \sum_{j_1,\dots,j_m} q_{j_1\dots j_m} x_1^{j_1}
\dots x_m^{j_m},
\]
then $v$ has entries $q_{j_1\dots j_m} X_1^{j_1} \dots X_m^{j_m}$, and
\begin{align*}
|Q(r_1, \dots, r_m)| &\leq \sum_{j_1,\dots,j_m} |q_{j_1\dots j_m}|
|r_1|^{j_1} \dots |r_m|^{j_m}\\
&\leq  \sum_{j_1,\dots,j_m} |q_{j_1\dots j_m}|
X_1^{j_1} \dots X_m^{j_m}\\
&= |v|_1.
\end{align*}
Thus, every vector $v \in L$
satisfying $|v|_1 < N^{\beta k}$ gives a polynomial relation between
$r_1,\dots,r_m$.

It is straightforward to compute the dimension and determinant of the
lattice:
\[
\dim L = \binom{t+m}{m},
\]
and
\[
\det L = (X_1 \dots X_m)^{\binom{t+m}{m}\frac{t}{m+1}} N^{\binom{k+m}{m} \frac{k}{m+1}}.
\]
To compute the determinant, we can choose a monomial ordering so that
the basis matrix for this lattice is upper triangular; then the
determinant is simply the product of the terms on the diagonal.

Now we apply LLL lattice basis reduction to $L$.  Because all the vectors in
$L$ are integral, the $m$ shortest vectors $v_1,\dots,v_m$ in the
LLL-reduced basis satisfy
\[
|v_1| \leq \dots \leq |v_m| \leq 2^{(\dim L )/4} (\det L)^{1/(\dim L+1-m)}
\]
(see Theorem 2 in \cite{herrmann2008solving}), and $|v|_1 \le
\sqrt{\dim L} \,|v|$ by Cauchy-Schwarz, so we know that the
corresponding polynomials $Q$ satisfy
\[
|Q(r_1, \dots, r_m)| \leq \sqrt{\dim L}\, 2^{(\dim L)/4} (\det L)^{1/(\dim L+1-m)}.
\]
If
\begin{equation}
\label{eq:desired}
\sqrt{\dim L}\, 2^{(\dim L)/4} \det L^{1/(\dim L+1-m)} < N^{\beta k},
\end{equation}
then we can conclude that $Q(r_1,\dots,r_m)=0$.

If $t$ and $k$ are large, then we can approximate $\binom{t+m}{m}$ with
$t^m/m!$ and $\binom{k+m}{m}$ with $k^m/m!$.  The $\sqrt{\dim L}$
factor plays no significant role asymptotically, so we simply omit it
(the omission is not difficult to justify). After taking a logarithm
and simplifying slightly, our desired equation \eqref{eq:desired}
becomes
\[
\frac{t^m}{4km!} + \frac{1}{1-\frac{(m-1)m!}{t^m}}\left(\frac{m\log_2 X}{m+1} \frac{t}{k} +
\frac{\log_2 N}{m+1} \frac{k^m}{t^m}\right) < \beta \log_2 N, 
\]
where $X$ denotes the geometric mean of $X_1,\dots,X_m$.

The $t^m/(4km!)$ and $(m-1)m!/t^m$ terms are nuisance factors, and once
we optimize the parameters they will tend to zero asymptotically. We
will take $t \approx \beta^{-1/m} k$ and $\log X \approx
\beta^{(m+1)/m} \log N$.  Then
\[
\frac{m\log X}{m+1} \frac{t}{k} +
\frac{\log N}{m+1} \frac{k^m}{t^m} \approx
\frac{m}{m+1} \beta \log N + \frac{1}{m+1} \beta \log N
= \beta \log N.
\]
By setting $\log X$ slightly less than this bound (by a $1+o(1)$
factor), we can achieve the desired inequality, assuming that the
$1-(m-1)!/t^m$ and $t^m/(4 k m!)$ terms do not interfere.  To ensure
that they do not, we take $t \gg m$ and $t^m \ll \beta \log N$ as $N
\to \infty$.  Note that then $\dim L \le \beta \log N$, which is
bounded independently of $m$.

Specifically, when $N$ is large we can take
\[
t = \left\lfloor\frac{(\beta \log N)^{1/m}}{(\beta^2 \log N)^{1/(2m)}}\right\rfloor
\]
and
\[
k = \lfloor \beta^{1/m} t \rfloor \approx (\beta^2 \log N)^{1/(2m)}.
\]
With these parameter settings, $t$ and $k$ both tend to infinity as $N
\to \infty$, because $\beta^2 \log N \to \infty$, and they satisfy the
necessary constraints. We do not recommend using these parameter
settings in practice; instead, one should choose $t$ and $k$ more
carefully. However, these choices work asymptotically.  Notice that
with this approach, $\beta^2 \log N$ must be large enough to allow
$t/k$ to approximate $\beta^{-1/m}$.  This is a fundamental issue, and
we discuss it in more detail in the next subsection.

The final step of the proof is to solve the system of equations defined
by the $m$ shortest vectors in the reduced basis to learn
$r_1,\dots,r_m$. One way to do this is to repeatedly use resultants to
eliminate variables; alternatively, we can use Gr\"{o}bner bases. See,
for example, Chapter~3 of \cite{IVA}.

One obstacle is that the equations may be not algebraically
independent, in which case we will not have enough information to
complete the solution.  In the experiments summarized in
Section~\ref{section:implementation}, we sometimes encountered cases
when the $m$ shortest vectors were algebraically dependent. However, in
every case the vectors represented either (1) irreducible,
algebraically independent polynomials, or (2) algebraically dependent
polynomials that factored easily into polynomials which all had the
desired properties.  Thus when the assumption of algebraic dependence
failed, it failed because there were fewer than $m$ independent factors
among the $m$ shortest relations.  In these cases, there were always
more than $m$ vectors of $\ell_1$ norm less than $N^{\beta k}$, and we
were able to complete the solution by using all these vectors.  This
behavior appears to depend sensitively on the optimization of the
parameters $t$ and $k$.

\subsection{The $\beta^2 \log N \gg 1$ requirement}

The condition that $\beta^2 \log N \gg 1$ is not merely a convenient
assumption for the analysis.  Instead, it is a necessary hypothesis for
our approach to work at all when using a lattice basis reduction
algorithm with an exponential approximation factor. In previous papers
on these lattice-based techniques, such as \cite{coppersmith1997small}
or \cite{howgrave2001approximate}, this issue seemingly does not arise,
but that is because it is hidden in a degenerate case.  When $m=1$, we
are merely ruling out the cases when the bound $N^{\beta^2}$ on the
perturbations is itself bounded, and in those cases the problem can be
solved by brute force.

To see why a lower bound on $\beta^2 \log N$ is necessary, we can start
with \eqref{eq:desired}.  For that equation to hold, we must at least
have $2^{(\dim L)/4} < N^{\beta k}$ and $(\det L)^{1/(\dim L)} <
N^{\beta k}$, and these inequalities imply that
\[
\frac{1}{4} \binom{t+m}{m} < \beta k \log_2 N
\]
and
\[
\frac{\binom{k+m}{m} \log_2 N}{\binom{t+m}{m} (m+1)} < \beta \log_2 N.
\]
Combining them with $\binom{k+m}{m} > k$ yields
\[
\frac{1}{4 (m+1)} < \beta^2 \log_2 N,
\]
so we have an absolute lower bound for $\beta^2 \log N$.  Furthermore,
one can check that in order for the $2^{(\dim L)/4}$ factor to become
negligible compared with $N^{\beta k}$, we must have $\beta^2 \log N
\gg 1$.

Given a lattice basis reduction algorithm with approximation factor
$2^{(\dim L)^{\varepsilon}}$, we could replace $t^m$ with
$t^{\varepsilon m}$ in the nuisance term coming from the approximation
factor.  Then the condition $t^m \ll \beta \log N$ would become
$t^{\varepsilon m} \ll \beta \log N$, and if we combine this with $k
\approx \beta^{1/m} t$, we find that
\[
k^{\varepsilon m} \approx \beta^\varepsilon t^{\varepsilon m}
\ll \beta^{1+\varepsilon} \log N.
\]
Because $k \ge 1$, the condition $\beta^{1+\varepsilon} \log N \gg 1$
is needed, and then we can take
\[
t = \left\lfloor\frac{(\beta \log N)^{1/(\varepsilon m)}}
{(\beta^{1+\varepsilon} \log N)^{1/(2\varepsilon m)}}\right\rfloor
\]
and
\[
k = \lfloor \beta^{1/m} t \rfloor \approx
(\beta^{1+\varepsilon} \log N)^{1/(2\varepsilon m)}.
\]

\subsection{Theorem~\ref{pseudo:gacd}}
\label{section:gacd}

The algorithm for Theorem~\ref{pseudo:gacd} is identical to the above,
except that we do not have an exact $N$, so we omit all vectors
involving $N$ from the construction of the lattice $L$.

The matrix of coefficients is no longer square, so we have to do more
work to bound the determinant of the lattice.
Howgrave-Graham~\cite{howgrave2001approximate} observed in the
two-variable case that the determinant is preserved even under
non-integral row operations, and he used a non-integral transformation
to hand-reduce the matrix before bounding the determinant as the
product of the $\ell_2$ norms of the basis vectors; furthermore, the
$\ell_2$ norms are bounded by $\sqrt{\dim L}$ times the $\ell_\infty$
norms.

The non-integral transformation that he uses is based on the relation
\[
(x_i-a_i) - \frac{a_i}{a_1} (x_1-a_1) = x_i - \frac{a_i}{a_1} x_1.
\]
By adding a multiple of $f(x)(x_1-a_1)$, one can reduce $f(x)(x_i-a_i)$
to $f(x)( x_i - \frac{a_i}{a_1} x_1)$.  The advantage of this is that
if $x_1 \approx x_i$ and $a_1 \approx a_i$, then $x_i - \frac{a_i}{a_1}
x_1$ may be much smaller than $x_i-a_i$ was.  The calculations are
somewhat cumbersome, and we will omit the details (see
\cite{howgrave2001approximate} for more information).

When $a_1,\dots,a_m$ are all roughly $N$ (as in
Theorem~\ref{pseudo:gacd}), we get the following values for the
determinant and dimension in the $m$-variable case:
\[
\det L \leq (N/X)^{\binom{k+m-1}{m}(t-k+1)}
X^{m\left(\binom{t+m}{m}\frac{t}{m+1} -
  \binom{k-1+m}{m}\frac{k-1}{m+1}\right)}
\]
and
\[
\dim L = \binom{t+m}{m} - \binom{k-1+m}{m}.
\]
To optimize the resulting bound, we take $t \approx (m/\beta)^{1/(m-1)}
k$.

\section{Applications to fully homomorphic encryption}
\label{section:fail}

In \cite{vandijkgentryhalevivaikuntanathan2010}, the authors build a
fully homomorphic encryption system whose security relies on several
assumptions, among them the hardness of computing an approximate common
divisor of many integers.  This assumption is used to build a simple
``somewhat homomorphic'' scheme, which is then transformed into a fully
homomorphic system under additional hardness assumptions.  In this
section, we use our algorithm for computing approximate common divisors
to provide a more precise understanding of the security assumption
underlying this somewhat homomorphic scheme, as well as the related
cryptosystem of \cite{coron2011fully}.

For ease of comparison, we will use the notation from the above two
papers (see Section~3 of \cite{vandijkgentryhalevivaikuntanathan2010}).
Let $\gamma$ be the bit length of $N$, $\eta$ be the bit length of $p$,
and $\rho$ be the bit length of each $r_i$. Using our algorithm, we can
find $r_1,\dots,r_m$ and the secret key $p$ when
\[
\rho \leq \gamma \beta^{(m+1)/{m}}.
\]
Substituting in $\beta = \eta/\gamma$, we obtain
\[
\rho^m \gamma \leq \eta^{m+1}.
\]
The authors of \cite{vandijkgentryhalevivaikuntanathan2010} suggest as
a ``convenient parameter set to keep in mind'' to set $\rho = \lambda$,
$\eta = \lambda^2$, and $\gamma = \lambda^5$.  Using $m>3$ we would be
able to solve this parameter set, if we did not have the barrier that
$\eta^2$ must be much greater than $\gamma$.

As pointed out in Section~\ref{subsub:acd}, this barrier would no
longer apply if we could improve the approximation factor for lattice
basis reduction.  If we could improve the approximation factor to
$2^{(\dim L)^\varepsilon}$, then the barrier would amount to
$\beta^{1+\varepsilon} \lambda^5 \gg 1$, where $\beta = \eta/\gamma =
\lambda^{-3}$.  If $\varepsilon < 2/3$, then this would no longer be an
obstacle.  Given a $2^{(\dim L)^{2/3}/\log \dim L}$ approximation
factor, we could take $m=4$, $k=1$, and $t = \lfloor 3 \lambda^{3/4}
\rfloor$ in the notation of Section~\ref{section:acd}. Then
\eqref{eq:desired} holds, and thus the algorithm works, for all
$\lambda \ge 300$.

One might try to achieve these subexponential approximation factors by
using blockwise lattice reduction techniques~\cite{gama2008finding}.
For an $n$-dimensional lattice, one can obtain an approximation factor
of roughly $\kappa^{n/\kappa}$ in time exponential in $\kappa$. For the
above parameter settings, the lattice will have dimension on the order
of $\lambda^3$, and even a $2^{n^{2/3}}$ approximation will require
$\kappa > n^{1/3} = \lambda$, for a running time that remains
exponential in $\lambda$. (Note that for these parameters, using a
subexponential-time factoring algorithm to factor the modulus in the
``partial'' approximate common divisor problem is super-exponential in
the security parameter.)

In general, if we could achieve an approximation factor of $2^{(\dim
L)^{\varepsilon}}$ for arbitrarily small $\varepsilon$, then we could
solve the approximate common divisor problem for parameters given by
any polynomials in $\lambda$.  Furthermore, as we will see in
Section~\ref{section:implementation}, the LLL algorithm performs better
in practice on these problems than the theoretical analysis suggests.

\section{Multi-polynomial reconstruction}
\label{section:parvaresh-vardy}

\subsection{Polynomial lattices}

For Theorem~\ref{multi-polynomial} and
Theorem~\ref{multi-polynomial-rigorous}, we can use almost exactly the
same technique, but with lattices over the polynomial ring $F[z]$
instead of the integers.

By a $d$-dimensional lattice $L$ over $F[z]$, we mean the $F[z]$-span
of $d$ linearly independent vectors in $F[z]^d$.  The degree $\deg v$
of a vector $v$ in $L$ is the maximum degree of any of its components,
and the determinant $\det L$ is the determinant of a basis matrix
(which is well-defined, up to scalar multiplication).

The polynomial analogue of lattice basis reduction produces a basis
$b_1,\dots,b_d$ for $L$ such that
\[
\deg(b_1) + \cdots + \deg (b_d) = \deg \det L.
\]
Such a basis is called a reduced basis (sometimes column or
row-reduced, depending on how the vectors are written), and it can be
found in polynomial time; see, for example, Section~6.3 in
\cite{Kailath}. If we order the basis so that $\deg(b_1) \le \dots \le
\deg(b_d)$, then clearly
\[
\deg(b_1) \le \frac{\deg \det L}{d},
\]
and more generally
\[
\deg(b_i) \le \frac{\deg \det L}{d-(i-1)},
\]
because
\[
\deg\det L - (d-(i-1))\deg(b_i) = \sum_{j=1}^d \deg(b_j) - \sum_{j=i}^d \deg(b_i) \ge 0.
\]
These inequalities are the polynomial analogues of the vector length
bounds in LLL-reduced lattices, but notice that the exponential
approximation factor does not occur.  See \cite{ideal-coppersmith} for
more information about this analogy, and \cite{DGH} for applications
that demonstrate the superior performance of these methods in practice.

\subsection{Theorems~\ref{multi-polynomial} and~\ref{multi-polynomial-rigorous}}

In the polynomial setting, we will choose $Q(x_1,\dots,x_m)$ to be a
linear combination (with coefficients from $F[z]$) of the polynomials
\[
(x_1-f_1(z))^{i_1} \dots (x_m - f_m(z))^{i_m} N(z)^\ell
\]
with $i_1 + \dots + i_m \le t$ and $\ell = \max(k - \sum_j i_j,0)$. We
define the lattice $L$ to be spanned by the coefficient vectors of
these polynomials, but with $x_i$ replaced with $z^{\ell_i} x_i$ to
incorporate the bound on $\deg g_i$, much as we replaced $x_i$ with
$X_ix_i$ in Section~\ref{section:acd}.

As before, we can easily compute the dimension and determinant of $L$:
\[
\dim L = \binom{t+m}{m}
\]
and
\[
\deg \det L = (\ell_1 + \dots + \ell_m) \binom{t+m}{m}\frac{t}{m+1} + n
\binom{k+m}{m} \frac{k}{m+1},
\]
where $n = \deg N(z)$.

Given a polynomial $Q(x_1,\dots,x_m)$ corresponding to a vector $v \in
L$, we can bound $\deg Q(g_1(z),\dots,g_m(z))$ by $\deg v$.
Specifically, suppose
\[
Q(x_1,\dots,x_m) = \sum_{j_1,\dots,j_m} q_{j_1\dots j_m}(z) x_1^{j_1}
\dots x_m^{j_m};
\]
then $v$ is the vector whose entries are $q_{j_1\dots j_m}(z) z^{j_1\ell_1+\dots+j_m\ell_m}$,
and
\begin{align*}
\deg Q(g_1(z), \dots, g_m(z))
&\leq \max_{j_1,\dots,j_m} (\deg q_{j_1\dots j_m}(z) + j_1 \deg g_1(z) + \dots +j_m \deg g_m(z))\\
&\leq \max_{j_1,\dots,j_m} (\deg q_{j_1\dots j_m}(z) + j_1 \ell_1 + \dots +j_m \ell_m)\\
&= \deg v.
\end{align*}

Let $v_1,\dots,v_{\dim L}$ be a reduced basis of $L$, arranged in
increasing order by degree. If
\begin{equation} \label{poly:ineq}
\frac{\deg \det L}{\dim L - (m-1)} < \beta k n,
\end{equation}
then each of $v_1,\dots,v_m$ yields a polynomial relation $Q_i$ such that
\[
Q_i(g_1(z),\dots,g_m(z))=0,
\]
because by the construction of the lattice,
$Q_i(g_1(z),\dots,g_m(z))$ is divisible by the $k$-th power of an
approximate common divisor of degree $\beta n$, while
\[
\deg Q_i(g_1(z),\dots,g_m(z)) \le \deg v_i < \beta k n.
\]
Thus, we must determine how large
$\ell_1+\dots+\ell_m$ can be, subject to the inequality \eqref{poly:ineq}.

If we set $t \approx k \beta^{-1/m}$ and
\[
\frac{\ell_1 + \dots + \ell_m}{m} < n \beta^{(m+1)/m},
\]
then inequality~\eqref{poly:ineq} is satisfied when $t$ and $k$ are
sufficiently large.  Because there is no analogue of the LLL
approximation factor in this setting, we do not have to worry about $t$
and $k$ becoming too large (except for the obvious restriction that
$\dim L$ must remain polynomially bounded), and there is no lower bound
on $\beta$.  Furthermore, we require no $1+o(1)$ factors, because all
degrees are integers and all the quantities we care about are rational
numbers with bounded numerators and denominators; thus, any
sufficiently close approximation might as well be exact, and we can
achieve this when $t$ and $k$ are polynomially large.

\section{Higher degree polynomials}

It is possible to generalize the results in the previous sections to
find solutions of a system of higher degree polynomials modulo divisors
of $N$.

\begin{theorem} 
\label{pseudo:highdeg} Given a positive integer $N$ and $m$ monic
polynomials $h_1(x),\dots,h_m(x)$ over the integers, of degrees
$d_1,\dots,d_m$, and given any $\beta \gg 1/\sqrt{\log N}$ and bounds
$X_1,\dots,X_m$, we can find all $r_1,\dots,r_m$ such that
\[
\gcd(N, h_1(r_1), \dots, h_m(r_m)) \ge N^\beta
\]
and $|r_i| \le X_i$, provided that
\[
\sqrt[m]{X_1^{d_1} \dots X_m^{d_m}} < N^{(1+o(1))\beta^{(m+1)/m}}
\]
and that the algebraic independence hypothesis holds. The algorithm
runs in polynomial time for fixed $m$.
\end{theorem}

The $m=1$ case does not require the algebraic independence hypothesis,
and it encompasses both Howgrave-Graham and Coppersmith's theorems
\cite{howgrave2001approximate,coppersmith1997small}; it first appeared
in \cite{may2003new}.

When $X_1 = \dots = X_m$, the bound becomes
$N^{\beta^{(m+1)/m}/\bar{d}}$, where $\bar{d} = (d_1+\dots+d_m)/m$ is
the average degree.

\begin{theorem} 
\label{pseudo:highdeg-poly} Given a polynomial $N(z)$ and $m$ monic
polynomials $h_1(x),\dots,h_m(x)$ over $F[z]$, of degrees
$d_1,\dots,d_m$ in $x$, and given degree bounds $\ell_1,\dots,\ell_m$,
we can find all $g_1(z),\dots,g_m(z)$ in $F[z]$ such that
\[
\deg \gcd(N(z), h_1(g_1(z)), \dots, h_m(g_m(z))) \ge \beta \deg N(z)
\]
and $\deg g_i(z) \le \ell_i$, provided that
\[
\frac{\ell_1d_1 + \dots + \ell_md_m}{m} < \beta^{(m+1)/m} \deg N(z)
\]
and that the algebraic independence hypothesis holds. The algorithm
runs in polynomial time for fixed $m$.
\end{theorem}

The algorithms are exactly analogous to those for the degree $1$ cases,
except that $x_i-a_i$ (or $x_i - f_i(z)$) is replaced with $h_i(x_i)$.

\section{Implementation}
\label{section:implementation}

We implemented the number-theoretic version of the partial approximate
common divisor algorithm using Sage \cite{sage}.  We used Magma
\cite{Magma} to do the LLL and Gr\"obner basis calculations.

\begin{table}[tb]
\centering
\begin{tabular}{cccccccccc}
\toprule
$m$ & $\log_2 N$ & $\log_2 p$ & $\log_2 r$ & $t$ & $k$ & $\dim L$ & LLL &
Gr\"{o}bner & LLL factor \\
\midrule
1 & 1000 & 200 & 36 & 41 & 8 & 42 & 12.10s & -- & 1.037 \\
\textit{1} & \textit{1000} & \textit{200} & \textit{39} & \textit{190} & \textit{38} & \textit{191} \\
1 & 1000 & 400 & 154 & 40 & 16 & 41 & 34.60s & -- & 1.023 \\
1 & 1000 & 400 & 156 & 82 & 33 & 83 & 4554.49s & --  & 1.029 \\
\textit{1} & \textit{1000} & \textit{400} & \textit{159} & \textit{280}
& \textit{112} & \textit{281} \\
\midrule
2 & 1000 & 200 & 72 & 9    & 4 & 55  & 25.22s & 0.94s & 1.030 \\
\textit{2} & \textit{1000} & \textit{200} & \textit{85} & \textit{36}
& \textit{16} & \textit{703} \\
2 & 1000 & 400 & 232 & 10 & 6 & 66 & 126.27s & 5.95s & 1.038 \\
 2 & 1000 & 400 & 238 & 15 & 9 & 136 & 15720.95s & 25.86s & 1.019 \\
\textit{2} & \textit{1000} & \textit{400} & \textit{246} & \textit{46}
& \textit{29} & \textit{1128} \\
\midrule
3 & 1000 & 200 & 87   & 5 & 3 & 56 & 18.57s & 1.20s &  1.038 \\
\textit{3} & \textit{1000} & \textit{200} & \textit{102} & \textit{14}
& \textit{8} & \textit{680} \\
3 & 1000 & 400 & 255 & 4 & 3 & 35 & 2.86s  &      2.13ss & 1.032 \\
3 & 1000 & 400 & 268 & 7 & 5 & 120 & 1770.04s &  25.43s &  1.040 \\
\textit{3} & \textit{1000} & \textit{400} & \textit{281} & \textit{19}
& \textit{14} & \textit{1540} \\
\midrule
4 & 1000 & 200 & 94 & 3 & 2 & 35 & 1.35s &  0.54s & 1.028 \\
\textit{4} & \textit{1000} & \textit{200} & \textit{111} & \textit{8}
& \textit{5} & \textit{495} \\
4 & 1000 & 400 & 279 & 4 & 3 & 70 & 38.32s & 9.33s & 1.035\\
\textit{4} & \textit{1000} & \textit{400} & \textit{293} & \textit{10}
& \textit{8} & \textit{1001} \\
\midrule
5 & 1000 & 200 & 108 & 3 & 2 & 56 & 7.35s & 1.42s &  1.035 \\
5 & 1000 & 200 & 110 & 4 & 3 & 126 & 738.57s & 7.28s & 1.037 \\
5 & 1000 & 400 & 278 & 3 & 2 & 56 & 1.86s & 0.90s* & 0.743 \\
\midrule
6 & 1000 & 200 & 115 & 3 & 2 & 84 & 31.51s & 3.16s & 1.038 \\
6 & 1000 & 400 & 297 & 3 & 2 & 84 & 3.97s & 1.34s* & 0.586 \\
\midrule
7 & 1000 & 200 & 120 & 3 & 2 & 120 & 203.03s & 7.73s & 1.046 \\
7 & 1000 & 400 & 311 & 3 & 2 & 120 & 12.99s & 2.23s* & 0.568 \\
\midrule
12 & 1000 & 400 & 347 & 1 & 1 & 13 & 0.01s & 0.52s &
$1.013$ \\
18 & 1000 & 400 & 364 & 1 & 1 & 19   & 0.03s &     1.08s &
$1.032$ \\
24 & 1000 & 400 & 372 & 1 & 1 & 25   & 0.04s &     1.93s&
$1.024$ \\
48 & 1000 & 400 & 383 & 1 & 1 & 49 & 0.28s & 8.37s & $1.030$ \\
96 & 1000 & 400 & 387 & 1 & 1 & 97 & 1.71s & 27.94s & $1.040$ \\
\bottomrule
\end{tabular}
\medskip
\caption{Experimental results from our implementation of the integer
  partial approximate common divisor algorithm, with sample parameters for more
  extreme calculations in italics.}
\label{table:runningtimes}
\end{table}

We solved the systems of equations by computing a Gr\"obner basis with
respect to the lexicographic monomial ordering, to eliminate variables.
Computing a Gr\"obner basis can be extremely slow, both in theory and
in practice. We found that it was more efficient to solve the equations
modulo a large prime, to limit the bit length of the coefficients in
the intermediate and final results.  Because $r_1,\dots,r_m$ are
bounded in size, we can simply choose a prime larger than $2 \max_i
|r_i|$.

We ran our experiments on a computer with a 3.30 GHz quad-core Intel
Core i5 processor and 8 GB of RAM.  Table~\ref{table:runningtimes}
shows a selection of sample running times for various parameter
settings.  For comparison, the table includes the $m=1$ case, which is
Howgrave-Graham's algorithm.  The italicized rows give example lattice
dimensions for larger inputs to illustrate the limiting behavior of the
algorithm.

The performance of the algorithm depends on the ratio of $t$ to
$k$, which should be approximately $\beta^{-1/m}$.  Incorrectly
optimized parameters often perform much worse than correctly
optimized parameters.  For example, when $m=3$, $\log_2 N =
1000$, and $\log_2 p = 200$, taking $(t,k) = (4,2)$ can handle
$84$-bit perturbations $r_i$, as one can see in
Table~\ref{table:runningtimes}, but taking $(t,k)=(4,3)$ cannot
even handle $60$ bits.

For large $m$, we experimented with using the non-optimized parameters
$(t,k)=(1,1)$, as reported in Table~\ref{table:runningtimes}.  For the
shortest vector only, the bounds would replace the exponent
$\beta^{(m+1)/m}$ with $(m+1)\beta/m - 1/m$, which is its tangent line
at $\beta=1$.  This bound is always worse, and it is trivial when
$\beta \le 1/(m+1)$, but it still approaches the optimal exponent
$\beta$ for large $m$. Our analysis does not yield a strong enough
bound for the $m$-th largest vector, but in our experiments the vectors
found by LLL are much shorter than predicted by the worst-case bounds,
as described below. Furthermore, the algorithm runs extremely quickly
with these parameters, because the lattices have lower dimensions and
the simultaneous equations are all linear.

The last column of the table, labeled ``LLL factor,'' describes
the approximation ratio obtained by LLL in the experiment.
Specifically, LLL factor $\lambda$ means
\[
|v_m| \approx \lambda^{\dim L} (\det L)^{1/(\dim L)},
\]
where $v_m$ is the $m$-th smallest vector in the LLL-reduced basis for
$L$.  Empirically, we find that all of the vectors in the reduced basis
are generally quite close in size, so this estimate is more appropriate
than using $1/(\dim L - (m-1))$ in the exponent (which we did in the
theoretical analysis, in order to get a rigorous bound). The typical
value is about $1.02$, which matches the behavior one would expect from
LLL on a randomly generated lattice, whose successive minima will all
be close to $\det L^{1/(\dim L)}$~\cite{nguyen2006lll}.  A handful of
our experimental parameters resulted in lattices whose shortest vectors
were much shorter than these bounds; this tended to correlate with a
small sublattice of algebraically dependent vectors.

Because of this, the reduced lattice bases in practice contain many
more than $m$ suitable polynomials, and we were able to speed up some
of the Gr\"obner basis calculations by including all of them in the
basis.  For example, the $m=7$, $\log_2 p=200$ Gr\"obner basis
calculation from Table~\ref{table:runningtimes} finished in 12 seconds
using 119 polynomials from the reduced lattice basis.

We marked cases where we encountered algebraically dependent relations
with an asterisk in Table~\ref{table:runningtimes}. In each case, we
were still able to solve the system of equations by including more
relations from the lattice (up to $\ell_1$ norm less than $N^{\beta
k}$) and solving this larger system.

\section*{Acknowledgements}

We thank Chris Umans and Alex Vardy for suggesting looking at
Parvaresh-Vardy codes, and Martin Albrecht for advice on computing
Gr\"obner bases in practice.  N.H.\ would like to thank MIT, CSAIL, and
Microsoft Research New England for their hospitality during the course
of this research.  This material is based upon work supported by an
AT\&T Labs Graduate Fellowship and by the National Science Foundation
under Award No.\ DMS-1103803.

\end{document}